%% file: manuscript.tex
\documentclass{jssc}

\def\textsubscript#1%
{$_{\text{#1}}$}

\def\cdd{\mbox{\boldmath$\cdot$}~}
\def\q{\quad} 

\newcommand{\rulex}{\hfill\rule{1mm}{3mm}}
\def\ay{\arraycolsep=1.5pt}
\def\d{\displaystyle}
\def\dfrac{\displaystyle\frac}
\input{amssym.def}
\include{graphix}

\input jsscN.tex
\begin{document}

%*************************************************************************************************************
% \biaoti{THE CAPITALIZED TITLE OF YOUR ARTICLE$^*$}{The list of authors' names with the LAST NAME capitalized
% and the authors' names should be separated by "\cdd"}{the first author's name \\ the first author's affiliation
% and Email address\\ the second author's name\\ the second author's affiliation. More can be listed like this.}
% {$^*$ The titles and numbers of the foundations that support this article.}
%*************************************************************************************************************
\biaoti{Discrete-time Mean-field Stochastic $H_2/H_\infty$ Control$^*$}%%%   Main Title of your paper  %%%
{\uppercase{zhang} Weihai \cdd \uppercase{ma}
Limin }%%% The names of the authors  %%%
{ZHANG Weihai\\
 College of Electrical Engineering and Automation, Shandong University of  Science and Technology,
 Qingdao 266590, Shandong Province,  China. \\
    Email:w\_hzhang@163.com. \\ %%% The address of the authors  %%%
MA Limin\\
 College of Electrical Engineering and Automation, Shandong University of  Science and Technology,
 Qingdao 266590, Shandong Province,  China. \\
 Science and Information College, Qingdao Agricultural University, Qingdao 266109, Shandong Province,  China.\\
    Email:hbchdmlm@163.com. }
{$^*$ This research was supported by NSF of China under Grant
No.61573227, the Research Fund for the Taishan Scholar Project of
Shandong Province of China.\\
{}}

%*************************************************************************************************************
%The submission date of your article. For example: \drd{Received: June 8, 2006}
%*************************************************************************************************************
%\drd{DOI: }{Received: x x 20xx}{ / Revised: x x 20xx}

%*************************************************************************************************************
% The page header of the article.
% \dshm{Year}{Volume}{The capitalized RUNNING HEAD of your article with less than 48 letters}{The capitalized
% AUTHORS list with $\cdot$ separating different names or one can type "The name of the first author et al."
% if there are more than 4 authors.}
%*************************************************************************************************************

%\dshm{20XX}{XX}{Discrete-time Mean-field Stochastic $H_2/H_\infty$ Control}{\uppercase{zhang} Weihai \cdd \uppercase{ma}
%Limin}

%*************************************************************************************************************
% \dab{The abstract}{Keywords}
%*************************************************************************************************************
%-------------------------------------------------------------------------
\Abstract{ The finite  horizon  $H_2/H_\infty$ control problem of
mean-field type for discrete-time systems is considered in this
paper. Firstly, we derive a mean-field stochastic bounded real lemma
(SBRL). Secondly, a sufficient condition for the solvability of
discrete-time  mean-field stochastic linear-quadratic (LQ) optimal
control  is presented. Thirdly, based on SBRL and LQ results, this
paper establishes a sufficient condition for the existence of
discrete-time stochastic $H_2/H_\infty$ control of mean-field type
via the solvability of
coupled matrix-valued equations. }      % the abstract

        % the keywords

%\MRSubClass{05B05, 05B25, 20B25}      % MR(2000) Subject Classification

%\baselineskip 15pt

\section{Introduction}

In this paper, we investigate a class of stochastic $H_2/H_\infty$
control problems of mean-field type such as
\begin{equation}
\label{E:1.1}
\begin{cases}
 \  x(k+1)= A(k)x(k)+\tilde{A}(k)\mathbf{E}x(k)+B(k)\nu(k)+\tilde{B}(k)\mathbf{E}\nu(k)&\\
\ \hspace{1.6cm}+[C(k)x(k)+\tilde{C}(k)\mathbf{E}x(k)+D(k)\nu(k)+\tilde{D}(k)\mathbf{E}\nu(k)]\omega(k)+F_{1}(k)u(k), & \\
 \ z(k)=\left[\begin{array}{c}
{\it \Phi}(k)x(k) \\
{\it \Psi}(k)u(k)
\end{array}\right], & \\
{\it \Psi}^{T}(k){\it \Psi}(k)=I,\q x(0)=x_{0}\in R^n. &
  \end{cases}
\end{equation}
Different from the classical stochastic $H_2/H_\infty$ control
problem, both the expectation $\mathbf{E}x(k)$ of the  system state
$x(k)$ and the expectation $\mathbf{E}\nu(k)$ of the disturbance
signal $\nu(k)$ appear in the state equation  (\ref{E:1.1}). Such an
equation is a discrete-time stochastic difference equation of
McKean-Vlasov type and is called  a mean-field stochastic difference
equation. In particular, the corresponding $H_2/H_\infty$ control
problem is referred to as a mean-field stochastic  $H_2/H_\infty$
control, which is a combination of stochastic $H_2/H_\infty$ control
problem and mean-field theory.

Mixed $H_2/H_\infty$ control has become one of the most popular
research issues in the last two decades, which has attracted
 considerable attention of many authors and has been widely applied to
 various fields;  see, e.g., \cite{ref1,ref2,ref3,ref4,ref5}  for the discussion  of deterministic systems.
 From  1998,  researchers have paid more attention to stochastic
 $H_2/H_\infty$ theory and made  great progress. We refer the reader to
\cite{ref6,ref7,ref8,ref9,ref10,ref11,ref12,ref13}  and the
references therein for details. The reference \cite{ref9} studied
 the discrete-time $H_\infty$ control with state and exogenous disturbance
 dependent noise,  while \cite{ref10} dealt with a class of discrete-time
  stochastic $H_2/H_\infty$ control with  additive disturbance.  Recently,
  the results of \cite{ref6} were extended to discrete-time systems in \cite{ref11,ref12}, where the  finite horizon
  and infinite horizon mixed $H_2/H_\infty$ control were investigated,  respectively. In particular,
   it was shown  that the solvability of
  the mixed $H_2/H_\infty$ control problem is equivalent to that of four coupled
  matrix-valued equations. In 2010, the results of \cite{ref11} were generalized to discrete-time stochastic
    systems with Markovian jumps and multiplicative noise in \cite{ref13}.

In recent years, mean-field theory  has attracted considerable
attention, which   is  developed to study the collective behaviors
resulting from individuals' mutual interactions in various physical
and sociological dynamical systems. In the survey paper \cite{ref22},
 three examples were presented to
use  mean-field approach to modelling in economics,  finance and
other related issues. Based on mean-field theory, mean-field term
presents  the interactions among elements, which approaches the
expected value when the number of agents goes to infinity. Similarly
to \cite{ref17}, suppose the dynamical equation  of particle
$i(i=1,...,M)$ is described as \ay
\begin{eqnarray}
x_{i}^{M}(k+1)&=&A(k)x_{i}^{M}(k)+\tilde{A}(k)\frac{1}{M}\d\sum_{j=1}^{M}x_{j}^{M}(k)+B(k)\nu(k)+F_{1}(k)u(k)
\nonumber\\
&&+[C(k)x_{i}^{M}(k)+\tilde{C}(k)\frac{1}{M}\d\sum_{j=1}^{M}x_{j}^{M}(k)+D(k)\nu(k)]\omega_{i}^{M}(k),
\label{eq vhvbh}
\end{eqnarray}
where  $\{\omega_{i}^{M}(k),k\in N\}$, i=1,...,M, are independent of
each other and have identical statistics law. Letting $M\rightarrow
\infty$, we obtain the following equation by the law of large
numbers:
\ay
\begin{eqnarray}
 \  x(k+1)&=&A(k)x(k)+\tilde{A}(k)\mathbf{E}x(k)+B(k)\nu(k)
 \nonumber\\
&&+[C(k)x(k)+\tilde{C}(k)\mathbf{E}x(k)+D(k)\nu(k)]\omega(k)+F_{1}(k)u(k),
\label {eq vvcfcgu}
 \end{eqnarray}
which is a special case of the state equation of (\ref{E:1.1}). The
continuous-time case of (\ref{E:1.1}) is a mean-field stochastic
differential equation (MFSDE), which is of great importance in
applications and was
 introduced as a stochastic toy model for the Vlasov kinetic equation
 of plasma in \cite{ref25}. Since about 1956,
 MFSDEs and their applications attracted many authors' attention;
 see \cite{ref18,ref19,ref20} for mean-field backward  stochastic
differential equations (BSDEs)  and  stochastic partial differential
equations (SPDEs),  \cite{ref14,ref15,ref16,ref21} for stochastic
maximal principle and  \cite{ref26,ref27} for LQG control of
mean-field type stochastic systems.   Specifically,
 continuous-time  and discrete-time  mean-field LQ problems  were  studied in \cite{ref24}
and  \cite{ref17}, respectively.   However, up to date, we know few
about $H_\infty$ or $H_2/H_\infty$ control results  for system
(\ref{E:1.1}). To this end, we will discuss the mean-field
stochastic $H_2/H_\infty$ control  of (\ref{E:1.1}) in this note.
Compared with the pure $H_\infty$ control, mixed $H_2/H_\infty$
control takes the robustness and optimality into account, and
therefore appears more attractive in practice \cite{ref30}.
Especially, infinite horizon $H_2/H_\infty$ control for
discrete-time time-varying Markov jump systems with multiplicative
noise has been applied to multiplier-accelerator macroeconomic
system\cite{ref31}.

In this paper, we will deal with the finite horizon stochastic
$H_2/H_\infty$ control, which extends the results  of \cite{ref11}
to discrete-time time-varying mean-field systems with $(x,
\nu)$-dependent noise. In consideration of the appearances of
$\mathbf{E}x(k)$ and $\mathbf{E}\nu(k)$ in system dynamics, we are
not able to solve mean-field stochastic $H_2/H_\infty$ control by
the same methods used in classical stochastic $H_2/H_\infty$
control.  Hence, this paper is by no means a trivial extension of
\cite{ref11}. In virtue of the representations $x(k)-\mathbf{E}x(k),
\mathbf{E}x(k);\nu(k)-\mathbf{E}\nu(k), \mathbf{E}\nu(k)$, we may
derive a SBRL for a class of discrete-time mean-field time-varying
systems with $(x, \nu)$-dependent noise. Roughly speaking, the SBRL
is a fundamental tool to handle the $H_\infty$ control and
estimation problems for stochastic systems. A sufficient condition
for the existence of discrete-time stochastic $H_2/H_\infty$ control
of mean-field type via the solvability of coupled matrix-valued
equations is provided as our main results. In this paper, the system
coefficients are time-varying, so the corresponding results of
time-invariant systems  are our special cases. In addition, a
recursive algorithm is provided to solve the coupled matrix-valued
equations.

The contribution of this paper is as follows: Section 2 gives a
mean-field  SBRL. A sufficient condition for the solvability of
discrete-time  mean-field  stochastic  LQ optimal control problem is
established in Section 3. Section 4 contains our main theorems.
 A recursive algorithm is provided to solve the coupled matrix-valued equations
accurately in Section 5. Finally, we end this paper in Section 6
with a brief conclusion.

For convenience, throughout the paper, we adopt the following
notations: $X^T$: the transpose of the  matrix $X$ or vector $X$.
$X\geq 0$ ($X>0$): $X$ is  positive semi-definite (positive
definite) symmetric matrix.   $R^{m}$: the $m$-dimensional real
vector space with the usual inner product.  $R^{m\times n}$: all
$m\times n$-dimensional matrices space with entries in $R$.
 $N_{K}=\{0,1,2,\cdots,K\}$.
 $N=\{0,1,2,\cdots,\}$.
 $H_{n}(R)$: the set of all real symmetric matrices.

\section{Stochastic Bounded Real Lemma}
In this section, our main purpose is to obtain a mean-field SBRL,
which is the footstone    in the study of stochastic $H_\infty$
control and estimation. Consider the following discrete-time
stochastic difference equation with $k\in N_{K}$:
\begin{equation}
\label{E:2.1}
\begin{cases}
 \  x(k+1)= A(k)x(k)+\tilde{A}(k)\mathbf{E}x(k)+B(k)\nu(k)+\tilde{B}(k)\mathbf{E}\nu(k)&\\
\ \hspace{1.6cm}+[C(k)x(k)+\tilde{C}(k)\mathbf{E}x(k)+D(k)\nu(k)+\tilde{D}(k)\mathbf{E}\nu(k)]\omega(k), & \\
 \ z(k)={\it \Phi}(k)x(k), &\\
 \ x(0)=x_{0}\in R^n, &
  \end{cases}
\end{equation}
where $x(k)\in R^n$, $\nu(k)\in R^l$, and $z(k)\in R^m$ are
respectively the system state, disturbance signal and controlled
output. $A(k),\tilde{A}(k),C(k),\tilde{C}(k)\in R^{n\times n}$,
$B(k),\tilde{B}(k),D(k),\tilde{D}(k)\in R^{n\times l}$, and ${\it \Phi}(k)\in
R^{m\times n}$ are given matrix-valued  functions. The initial value $x_{0}$ is assumed to be a deterministic vector.
$\mathbf{E}$ is the
expectation operator. $\{\omega(k),k\in N_{K}\}$ is a sequence of real random variables defined on a
complete probability space $\{{\it \Omega},\mathcal{F},\mu\}$, which is a wide sense stationary, second
order process with $\mathbf{E}(\omega(s))=0$ and
$\mathbf{E}(\omega(s)\omega(t))=\delta_{st}$, where $\delta_{st}$ is a
Kronecker function. Suppose $\omega(k)$ and $\nu(k)$ are
uncorrelated. Denote $\mathcal{F}_{k}$ the $\sigma$-algebra generated by
$\{\omega(t), t=0,1,\cdots,k\}$. Let $L^{2}({\it \Omega},R^{p})$ be the
space of $R^{p}$-valued square integrable random vectors, and
$l^{2}_{\omega}(N_{K},R^{p})$ denotes the space of all finite
sequences $y(k)\in L^{2}({\it \Omega},R^{p})$ that are $\mathcal{F}_{k-1}$
measurable for $k\in N_{K}$. The $l^{2}$ norm of
$l^{2}_{\omega}(N_{K},R^{p})$ is defined as
$$\|y(.)\|_{l^{2}_{\omega}(N_{K},R^{p})}=\bigg(\d\sum\limits_{k=0 }^{K} \mathbf{E}\|y(k)\|^{2}\bigg)^{1/2}.$$
For any $K\in N$ and $(x_{0},\nu(k))\in R^{n}\times
l^{2}_{\omega}(N_{K},R^{l})$, the unique solution of (\ref{E:2.1})
with initial value $x_{0}$ is described as $x(k;x_{0},\nu)$.

\begin{definition}\label{df:2.1}
The perturbed operator of system
(\ref{E:2.1}) is defined by
\[
L_{K}:l^{2}_{\omega}(N_{K},R^{l}) \rightarrow
l^{2}_{\omega}(N_{K},R^{m}),
\]
\[L_{K}(\nu(k)):={\it \Phi}(k)x(k;0,\nu),\q \forall \ \nu(k)\in l^{2}_{\omega}(N_{K},R^{l})\]
 with its norm
$$
\|L_{K}\|=\mathop {\sup }\limits_{\nu\in l^{2}_{\omega}(N_{K},R^{l}),\nu\neq 0}
\dfrac{\|z(k)\|_{l^{2}_{\omega}(N_{K},R^{m})}}{\|\nu(k)\|_{l^{2}_{\omega}(N_{K},R^{l})}}
=\mathop {\sup }\limits_{\nu\in l^{2}_{\omega}(N_{K},R^{l}),\nu\neq 0}
 \dfrac{\bigg(\d\sum\limits_{k=0 }^{K}\mathbf{E}\|{\it \Phi}(k)x(k;0,\nu)\|^{2}\bigg)^{1/2}}{\bigg(\d\sum\limits_{k=0
 }^{K}\mathbf{E}\|\nu(k)\|^{2}\bigg)^{1/2}}.
$$
\end{definition}

In this paper, we discuss  the  mean-field $H_2/H_\infty$ control,
which is a combination of mean-field theory and mixed $H_2/H_\infty$
control problem. In consideration of the appearance of
$\mathbf{E}x(k)$ and $\mathbf{E}\nu(k)$ in the system (\ref{E:2.1}),
we may solve mean-field stochastic $H_2/H_\infty$ control problem by
using the representation $x(k)-\mathbf{E}x(k),
\mathbf{E}x(k);\nu(k)-\mathbf{E}\nu(k), \mathbf{E}\nu(k)$, which is
different from classical $H_2/H_\infty$ control problem.
 Taking expectations in system (\ref{E:2.1}),
we have the system equations on  $\mathbf{E}x(k)$ and $x(k)-\mathbf{E}x(k)$ satisfying
\begin{equation}
\label{E:2.2}
\begin{cases}
 \  \mathbf{E}x(k+1)= \mathbb{A}(k)\mathbf{E}x(k)+\mathbb{B}(k)\mathbf{E}\nu(k),&\\
 \ \mathbf{E}x(0)=\mathbf{E}x_{0}\in R^{n}, &
  \end{cases}
\end{equation}
\begin{equation}
\label{E:2.3}
\begin{cases}
 \  x(k+1)-\mathbf{E}x(k+1)=[A(k)(x(k)-\mathbf{E}x(k))+B(k)(\nu(k)-\mathbf{E}\nu(k))]&\\
\ \hspace{3.6cm}+[C(k)(x(k)-\mathbf{E}x(k))+\mathbb{C}(k)\mathbf{E}x(k) & \\
 \ \hspace{3.6cm}+D(k)(\nu(k)-\mathbf{E}\nu(k))+\mathbb{D}(k)\mathbf{E}\nu(k)]\omega(k),& \\
  \ x(0)-\mathbf{E}x(0)=x_{0}-\mathbf{E}x_{0}=0.&
  \end{cases}
\end{equation}
Here and hereafter, $\mathbb{A}(k)=A(k)+\tilde{A}(k)$, $\mathbb{B}(k)=B(k)+\tilde{B}(k)$,
$\mathbb{C}(k)=C(k)+\tilde{C}(k)$, $\mathbb{D}(k)=D(k)+\tilde{D}(k)$.

Next, based on the above equations, we may arrive at a SBRL  step by
step.

\begin{lemma}\label{Le:2.2}
In system (\ref{E:2.1}), suppose $K\in N$ is given,
 $$P(0), \q  P(1), \q  P(2), \q \cdots,  \q P(N+1); \q  Q(0), \q Q(1), \q Q(2), \q \cdots, \q Q(N+1)$$
are arbitrary families of matrices in $H_{n}(R)$, then for any $x_{0}\in R^{n}$,
we have
\ay
\begin{eqnarray*}
&& \d\sum\limits_{k=0 }^{K} \mathbf{E} \left\{\left[\begin{array}{c}
x(k)-\mathbf{E}x(k) \\
\nu(k)-\mathbf{E}\nu(k)
\end{array}\right]^{T}M(P)
\left[\begin{array}{c}
x(k)-\mathbf{E}x(k) \\
\nu(k)-\mathbf{E}\nu(k)
\end{array}\right]\right\}
\nonumber\\
&&+\sum\limits_{k=0 }^{K}\left[\begin{array}{c}
\mathbf{E}x(k) \\
\mathbf{E}\nu(k)
\end{array}\right]^{T}S(P,Q)
\left[\begin{array}{c}
\mathbf{E}x(k) \\
\mathbf{E}\nu(k)
\end{array}\right]
\nonumber\\
=&&\mathbf{E}[(x(K+1)-\mathbf{E}x(K+1))^{T}P(K+1)(\cdots)]
\nonumber\\
&&+[\mathbf{E}x(K+1)]^{T}Q(K+1)[\mathbf{E}x(K+1)]-[\mathbf{E}x_{0}]^{T}Q(0)[\mathbf{E}x_{0}],
\end{eqnarray*}
where in the above and what follows, when we write $M^TR(\cdots)$ or
$M^TR[\cdots]$ for simplicity, we mean $(\cdots)=M$  or
$[\cdots]=M$. In addition,
{\small\ay
\begin{eqnarray*}
M(P)=
\left[\begin{array}{cc}
-P(k)+A(k)^{T}P(k+1)A(k) &  A(k)^{T}P(k+1)B(k)      \\
+C(k)^{T}P(k+1)C(k) &  +C(k)^{T}P(k+1)D(k)  \\
                    &                     \\
B(k)^{T}P(k+1)A(k) & B(k)^{T}P(k+1)B(k)  \\
+D(k)^{T}P(k+1)C(k) & +D(k)^{T}P(k+1)D(k)
\end{array}\right]
\end{eqnarray*}}
and{\small
\ay
\begin{eqnarray*}S(P,Q)=
\left[\begin{array}{cc}
-Q(k)+\mathbb{A}(k)^{T}Q(k+1)\mathbb{A}(k) & \mathbb{A}(k)^{T}Q(k+1)\mathbb{B}(k)\\
+\mathbb{C}(k)^{T}P(k+1)\mathbb{C}(k) & +\mathbb{C}(k)^{T}P(k+1)\mathbb{D}(k)   \\
                   &                     \\
\mathbb{B}(k)^{T}Q(k+1)\mathbb{A}(k) & \mathbb{B}(k)^{T}Q(k+1)\mathbb{B}(k) \\
+\mathbb{D}(k)^{T}P(k+1)\mathbb{C}(k) & +\mathbb{D}(k)^{T}P(k+1)\mathbb{D}(k)
\end{array}\right].
\end{eqnarray*}}
\end{lemma}

\proof Since $\omega(k)$ is independent of  $x_{0}$, $x(k)$
 and $v(k)$,
in view  of $\mathbf{E}\omega(k)=0$ and
$\mathbf{E}[\omega(k_{1})\omega(k_{2})]=\delta_{k_{1}k_{2}}$, we
have \ay
\begin{eqnarray*}
&& \mathbf{E}\{[A(k)(x(k)-\mathbf{E}x(k))+B(k)(\nu(k)-\mathbf{E}\nu(k))]^{T}P(k+1)
\nonumber\\
&&\times[C(k)(x(k)-\mathbf{E}x(k))+\mathbb{C}(k)\mathbf{E}x(k)
+D(k)(\nu(k)-\mathbf{E}\nu(k))+\mathbb{D}(k)\mathbf{E}\nu(k)]\omega(k)\}=0.
\end{eqnarray*}
So  equations  (\ref{E:2.2}) and (\ref{E:2.3}) lead to \ay
\begin{eqnarray}\label {E:2.4}
&&\mathbf{E}\left\{[x(k+1)-\mathbf{E}x(k+1)]^{T}P(k+1)[x(k+1)-\mathbf{E}x(k+1)]\right.\nonumber\\
&&-[x(k)-\mathbf{E}x(k)]^{T}P(k)[x(k)-\mathbf{E}x(k)]\}\nonumber\\
&=&\mathbf{E}\{[A(k)(x(k)-\mathbf{E}x(k))+B(k)(\nu(k)-\mathbf{E}\nu(k))]^{T}P(k+1)[\cdots]\nonumber\\
&&+[C(k)(x(k)-\mathbf{E}x(k))+\mathbb{C}(k)\mathbf{E}x(k)+D(k)(\nu(k)-\mathbf{E}\nu(k))+\mathbb{D}(k)\mathbf{E}\nu(k)]^{T}\nonumber\\
&&\times
P(k+1)[\cdots]-[x(k)-\mathbf{E}x(k)]^{T}P(k)[\cdots]\}\nonumber\\
&=&\mathbf{E}\{[x(k)-\mathbf{E}x(k)]^{T}[A(k)^{T}P(k+1)A(k)+C(k)^{T}P(k+1)C(k)-P(k)][\cdots]\nonumber\\
&&+[x(k)-\mathbf{E}x(k)]^{T}[A(k)^{T}P(k+1)B(k)+C(k)^{T}P(k+1)D(k)][\nu(k)-\mathbf{E}\nu(k)]\nonumber\\
&&+[\nu(k)-\mathbf{E}\nu(k)]^{T}[B(k)^{T}P(k+1)A(k)+D(k)^{T}P(k+1)C(k)][x(k)-\mathbf{E}x(k)]\nonumber\\
&&\left.+[\nu(k)-\mathbf{E}\nu(k)]^{T}[B(k)^{T}P(k+1)B(k)+D(k)^{T}P(k+1)D(k)][\cdots]\right\}\nonumber\\
&&+[\mathbf{E}x(k)]^{T}[\mathbb{C}(k)^{T}P(k+1)\mathbb{C}(k)][\mathbf{E}x(k)]
+[\mathbf{E}x(k)]^{T}[\mathbb{C}(k)^{T}P(k+1)\mathbb{D}(k)][\mathbf{E}\nu(k)]\nonumber\\
&&+[\mathbf{E}\nu(k)]^{T}[\mathbb{D}(k)^{T}P(k+1)\mathbb{C}(k)][\mathbf{E}x(k)]
+[\mathbf{E}\nu(k)]^{T}[\mathbb{D}(k)^{T}P(k+1)\mathbb{D}(k)][\mathbf{E}\nu(k)]
\end{eqnarray}
and
\ay
\begin{eqnarray}
&&\mathbf{E}\{[\mathbf{E}\lefteqn{ x(k+1)]^{T}Q(k+1)[\mathbf{E}x(k+1)]-[\mathbf{E}x(k)]^{T}Q(k)[\mathbf{E}x(k)]\}}\nonumber\\
&=&\mathbf{E}\{[\mathbb{A}(k)\mathbf{E}x(k)+\mathbb{B}(k)\mathbf{E}\nu(k)]^{T}\nonumber\\
&&\times
Q(k+1)[\mathbb{A}(k)\mathbf{E}x(k)+\mathbb{B}(k)\mathbf{E}\nu(k)]-[\mathbf{E}x(k)]^{T}Q(k)[\mathbf{E}x(k)]\}\nonumber\\
&=&[\mathbf{E}x(k)]^{T}[\mathbb{A}(k)^{T}Q(k+1)\mathbb{A}(k)-Q(k)][\mathbf{E}x(k)]
+[\mathbf{E}x(k)]^{T}[\mathbb{A}(k)^{T}Q(k+1)\mathbb{B}(k)][\mathbf{E}\nu(k)]\nonumber\\
&&+[\mathbf{E}\nu(k)]^{T}[\mathbb{B}(k)^{T}Q(k+1)\mathbb{A}(k)][\mathbf{E}x(k)]
+[\mathbf{E}\nu(k)]^{T}[\mathbb{B}(k)^{T}Q(k+1)\mathbb{B}(k)][\mathbf{E}\nu(k)].
\label {E:2.5}
\end{eqnarray}
Taking summation on both sides of  (\ref{E:2.4}) and (\ref{E:2.5})
over $k=0,1,2,\cdots,K$, respectively, we draw the conclusion
of this lemma.\rulex

\begin{lemma}\label{Le:2.3}
In system (\ref{E:2.1}), suppose $K\in N$ is given,
 $$P(0), P(1), P(2),\cdots, P(N+1); Q(0),Q(1),Q(2),\cdots,Q(N+1)$$
are arbitrary families of matrices in $H_{n}(R)$, then for any
$x_{0}\in R^{n}$, $\nu(k)\in l^{2}_{\omega}(N_{K},R^{l})$,
 we have
 \ay
\begin{eqnarray*}
J^{K}(x_{0},\nu)&=&\d\sum\limits_{k=0
}^{K}\mathbf{E}[\gamma^{2}\|\nu(k)\|^{2}-\|z(k)\|^{2}]
\nonumber\\
&=&\d\sum\limits_{k=0 }^{K} \mathbf{E} \left\{\left[\begin{array}{c}
x(k)-\mathbf{E}x(k) \\
\nu(k)-\mathbf{E}\nu(k)
\end{array}\right]^{T}\tilde{M}(P)
\left[\begin{array}{c}
x(k)-\mathbf{E}x(k) \\
\nu(k)-\mathbf{E}\nu(k)
\end{array}\right]\right\}
\nonumber\\
&&+\d\sum\limits_{k=0 }^{K} \left[\begin{array}{c}
\mathbf{E}x(k) \\
\mathbf{E}\nu(k)
\end{array}\right]^{T}\tilde{S}(P,Q)
\left[\begin{array}{c}
\mathbf{E}x(k) \\
\mathbf{E}\nu(k)
\end{array}\right]
\nonumber\\
&&-\mathbf{E}[(x(K+1)-\mathbf{E}x(K+1))^{T}P(K+1)(\cdots)]
\nonumber\\
&&+[\mathbf{E}x_{0}]^{T}Q(0)[\mathbf{E}x_{0}]-[\mathbf{E}x(K+1)]^{T}Q(K+1)[\mathbf{E}x(K+1)],
\end{eqnarray*}
where
\ay
\begin{eqnarray*}
\tilde{M}(P)=
\left[\begin{array}{cc}
-P(k)+A(k)^{T}P(k+1)A(k) &  A(k)^{T}P(k+1)B(k)      \\
+C(k)^{T}P(k+1)C(k) &  +C(k)^{T}P(k+1)D(k)  \\
-{\it \Phi}(k)^{T}{\it \Phi}(k) &                 \\
                   & \gamma^{2}I_{l}         \\
B(k)^{T}P(k+1)A(k) & +B(k)^{T}P(k+1)B(k)  \\
+D(k)^{T}P(k+1)C(k) & +D(k)^{T}P(k+1)D(k)
\end{array}\right]
\end{eqnarray*}
and
\ay
\begin{eqnarray*}
\tilde{S}(P,Q)=
\left[\begin{array}{cc}
-Q(k)+\mathbb{A}^{T}(k)Q(k+1)\mathbb{A}(k) & \mathbb{A}(k)^{T}Q(k+1)\mathbb{B}(k)\\
+\mathbb{C}(k)^{T}P(k+1)\mathbb{C}(k) & +\mathbb{C}(k)^{T}P(k+1)\mathbb{D}(k)   \\
-{\it \Phi}(k)^{T}{\it \Phi}(k) &                 \\
                   & \gamma^{2}I_{l}         \\
\mathbb{B}(k)^{T}Q(k+1)\mathbb{A}(k) & +\mathbb{B}(k)^{T}Q(k+1)\mathbb{B}(k) \\
+\mathbb{D}(k)^{T}P(k+1)\mathbb{C}(k) & +\mathbb{D}(k)^{T}P(k+1)\mathbb{D}(k)
\end{array}\right].
\end{eqnarray*}
\end{lemma}

\proof In view of $\nu(k)=\nu(k)-\mathbf{E}\nu(k)+\mathbf{E}\nu(k)$ and $x(k)=x(k)-\mathbf{E}x(k)+\mathbf{E}x(k)$, we have
\ay
\begin{eqnarray*}
\mathbf{E}[\gamma^{2}\|\nu(k)\|^{2}]&=&\lefteqn{\mathbf{E}\{[\nu(k)-\mathbf{E}\nu(k)+\mathbf{E}\nu(k)]^{T}\gamma^{2}[\nu(k)-\mathbf{E}\nu(k)+\mathbf{E}\nu(k)]\}}
\nonumber\\
&=&\mathbf{E}[(\nu(k)-\mathbf{E}\nu(k))^{T}\gamma^{2}I_{l}(\nu(k)-\mathbf{E}\nu(k))+(\mathbf{E}\nu(k))^{T}\gamma^{2}I_{l}(\mathbf{E}\nu(k))]
\end{eqnarray*}
and
\ay
\begin{eqnarray*}
\mathbf{E}[\|z(k)\|^{2}]&=&\lefteqn{
\mathbf{E}\{[{\it \Phi}(k)(x(k)-\mathbf{E}x(k)+\mathbf{E}x(k))]^{T}[{\it \Phi}(k)(x(k)-\mathbf{E}x(k)+\mathbf{E}x(k))]\}}\nonumber\\
&=&\mathbf{E}[(x(k)-\mathbf{E}x(k))^{T}{\it \Phi}(k)^{T}{\it \Phi}(k)(x(k)-\mathbf{E}x(k))+(\mathbf{E}x(k))^{T}{\it \Phi}(k)^{T}{\it \Phi}(k)(\mathbf{E}x(k))].
\end{eqnarray*}
From Lemma \ref{Le:2.2}, we have
\begin{eqnarray*}
J^{K}(x_{0},\nu)&=& \d\sum\limits_{k=0
}^{K}\mathbf{E}[\gamma^{2}\|\nu(k)\|^{2}-\|z(k)\|^{2}]\nonumber\\
&=&\d\sum\limits_{k=0
}^{K}\mathbf{E}[(\nu(k)-\mathbf{E}\nu(k))^{T}\gamma^{2}I_{l}(\nu(k)-\mathbf{E}\nu(k))+(\mathbf{E}\nu(k))^{T}\gamma^{2}I_{l}(\mathbf{E}\nu(k))
\nonumber\\
&&-(x(k)-\mathbf{E}x(k))^{T}{\it \Phi}(k)^{T}{\it \Phi}(k)(\cdots)-(\mathbf{E}x(k))^{T}{\it \Phi}(k)^{T}{\it \Phi}(k)(\mathbf{E}x(k))]\nonumber\\
&=&\d\sum\limits_{k=0 }^{K} \mathbf{E} \left\{\left[\begin{array}{c}
x(k)-\mathbf{E}x(k) \\
\nu(k)-\mathbf{E}\nu(k)
\end{array}\right]^{T}\tilde{M}(P)
\left[\begin{array}{c}
x(k)-\mathbf{E}x(k) \\
\nu(k)-\mathbf{E}\nu(k)
\end{array}\right]\right\}
\nonumber\\
&&+\sum\limits_{k=0 }^{K} \left[\begin{array}{c}
\mathbf{E}x(k) \\
\mathbf{E}\nu(k)
\end{array}\right]^{T}\tilde{S}(P,Q)
\left[\begin{array}{c}
\mathbf{E}x(k) \\
\mathbf{E}\nu(k)
\end{array}\right]
\nonumber\\
&&-\mathbf{E}[(x(K+1)-\mathbf{E}x(K+1))^{T}P(K+1)(\cdots)]
\nonumber\\
&&+[\mathbf{E}x_{0}]^{T}Q(0)[\mathbf{E}x_{0}]-[\mathbf{E}x(K+1)]^{T}Q(K+1)[\mathbf{E}x(K+1)],
\end{eqnarray*}
which ends the proof.\rulex

For convenience, we adopt the following notations:
\ay
\begin{eqnarray*}
&&L(P(k+1))= A(k)^{T}P(k+1)A(k)+C(k)^{T}P(k+1)C(k)-{\it \Phi}(k)^{T}{\it \Phi}(k),\nonumber \\
&&G(P(k+1))=A(k)^{T}P(k+1)B(k)+C(k)^{T}P(k+1)D(k), \nonumber \\
&&H(P(k+1))=\gamma^{2}I_{l}+B(k)^{T}P(k+1)B(k)+D(k)^{T}P(k+1)D(k),\nonumber\\
&&\tilde{L}(P(k+1),Q(k+1))=\mathbb{A}(k)^{T}Q(k+1)\mathbb{A}(k)
+\mathbb{C}(k)^{T}P(k+1)\mathbb{C}(k)-{\it \Phi}(k)^{T}{\it \Phi}(k),
\nonumber\\
&&\tilde{G}(P(k+1),Q(k+1))=\mathbb{A}(k)^{T}Q(k+1)\mathbb{B}(k)
+\mathbb{C}(k)^{T}P(k+1)\mathbb{D}(k),
\nonumber\\
&&\tilde{H}(P(k+1),Q(k+1))=\gamma^{2}I_{l}+\mathbb{B}(k)^{T}Q(k+1)\mathbb{B}(k)
+\mathbb{D}(k)^{T}P(k+1)\mathbb{D}(k).
\end{eqnarray*}

\begin{theorem} \label{th:2.4}
(SBRL) For mean-field type stochastic system (\ref{E:2.1}),
 we have $\|L_{K}\|<\gamma$ for some $\gamma>0$ and  $Q_{1}(k)\leq 0$ if the following constrained backward difference
equation
\begin{equation}
\label{E:2.6}
\begin{cases}
 \  P(k)=L(P(k+1))-G(P(k+1))H(P(k+1))^{-1}G(P(k+1))^{T}, &\\
 \  Q(k)=\tilde{L}(P(k+1),Q(k+1))-\tilde{G}(P(k+1),Q(k+1))&\\
 \ \hspace{1.3cm}\times\tilde{H}(P(k+1),Q(k+1))^{-1}\tilde{G}(P(k+1),Q(k+1))^{T},  & \\
 \ P(K+1)=Q(K+1)=0, &\\
 \ H(P(k+1))>0, \q \tilde{H}(P(k+1),Q(k+1))>0
  \end{cases}
\end{equation}
has a unique solution $(P_{1}(k),Q_{1}(k))$ .
\end{theorem}

\proof From Lemma \ref{Le:2.3} and $P(K+1)=Q(K+1)=0$, for $x_{0}=0$ we obtain
\ay
\begin{eqnarray*}
J^{K}(0,\nu)&=&
\d\sum\limits_{k=0 }^{K}\mathbf{E}[\gamma^{2}\|\nu(k)\|^{2}-\|z(k)\|^{2}]
\nonumber\\
&=&\d\sum\limits_{k=0 }^{K} \mathbf{E}
\left\{\left[\begin{array}{c}
x(k)-\mathbf{E}x(k) \\
\nu(k)-\mathbf{E}\nu(k)
\end{array}\right]^{T}\tilde{M}(P_{1})
\left[\begin{array}{c}
x(k)-\mathbf{E}x(k) \\
\nu(k)-\mathbf{E}\nu(k)
\end{array}\right]\right\}
\nonumber\\
&&+\d\sum\limits_{k=0 }^{K}
\left[\begin{array}{c}
\mathbf{E}x(k) \\
\mathbf{E}\nu(k)
\end{array}\right]^{T}\tilde{S}(P_{1},Q_{1})
\left[\begin{array}{c}
\mathbf{E}x(k) \\
\mathbf{E}\nu(k)
\end{array}\right].
\end{eqnarray*}
By completing squares method, we obtain for any $\nu(k)\in
l^{2}_{\omega}(N_{K},R^{l})$ with $\nu(k)\neq 0$, \ay
\begin{eqnarray*}
&&\d\sum\limits_{k=0 }^{K}\mathbf{E}\left\{\left[\begin{array}{c}
x(k)-\mathbf{E}x(k) \\
\nu(k)-\mathbf{E}\nu(k)
\end{array}\right]^{T}\tilde{M}(P_{1})
\left[\begin{array}{c}
x(k)-\mathbf{E}x(k) \\
\nu(k)-\mathbf{E}\nu(k)
\end{array}\right]\right\}
\nonumber\\
&=&\d\sum\limits_{k=0 }^{K}
\mathbf{E}\{(x(k)-\mathbf{E}x(k))^{T}[-P_{1}(k)+L(P_{1}(k+1))
\nonumber\\
&&-G(P_{1}(k+1))H(P_{1}(k+1))^{-1}G(P_{1}(k+1))^{T}](\cdots)\}
\nonumber\\
&&+\d\sum\limits_{k=0
}^{K}\mathbf{E}\{[(\nu(k)-\mathbf{E}\nu(k))-(\nu^{*}(k)-\mathbf{E}\nu^{*}(k))]^{T}
H(P_{1}(k+1))[\cdots]\}
\end{eqnarray*}
\ay
\begin{eqnarray*}
&=&\d\sum\limits_{k=0
}^{K}\mathbf{E}\{[(\nu(k)-\mathbf{E}\nu(k))-(\nu^{*}(k)-\mathbf{E}\nu^{*}(k))]^{T}H(P_{1}(k+1))[\cdots]\},
\end{eqnarray*}
where
$$\nu^{*}(k)-\mathbf{E}\nu^{*}(k)=-H(P_{1}(k+1))^{-1}G(P_{1}(k+1))^{T}(x(k)-\mathbf{E}x(k)).$$
In addition,
\ay
\begin{eqnarray*}
&&\sum\limits_{k=0 }^{K}
\left[\begin{array}{c}
\mathbf{E}x(k) \\
\mathbf{E}\nu(k)
\end{array}\right]^{T}\tilde{S}(P_{1},Q_{1})
\left[\begin{array}{c}
\mathbf{E}x(k) \\
\mathbf{E}\nu(k)
\end{array}\right]
\nonumber\\
&=&\d\sum\limits_{k=0 }^{K} \mathbf{E}x(k)^{T}[-Q_{1}(k)+\tilde{L}(P_{1}(k+1),Q_{1}(k+1))
\nonumber\\
&& -\tilde{G}(P_{1}(k+1),Q_{1}(k+1))\tilde{H}(P_{1}(k+1),Q_{1}(k+1))^{-1}\tilde{G}(P_{1}(k+1),Q_{1}(k+1))^{T}]\mathbf{E}x(k)
\nonumber\\
&& +\d\sum\limits_{k=0 }^{K}[\mathbf{E}\nu(k)+\tilde{H}(P_{1}(k+1),Q_{1}(k+1))^{-1}\tilde{G}(P_{1}(k+1),Q_{1}(k+1))^{T}\mathbf{E}x(k)]^{T}
\nonumber\\
&& \times \tilde{H}(P_{1}(k+1),Q_{1}(k+1))[\cdots]
\nonumber\\
&=&\d\sum\limits_{k=0
}^{K}[\mathbf{E}\nu(k)-\mathbf{E}\nu^{*}(k)]^{T}\tilde{H}(P_{1}(k+1),Q_{1}(k+1))[\cdots],
\end{eqnarray*}
where
$$\mathbf{E}\nu^{*}(k)=-\tilde{H}(P_{1}(k+1),Q_{1}(k+1))^{-1}\tilde{G}(P_{1}(k+1),Q_{1}(k+1))^{T}\mathbf{E}x(k).$$
So we have $J^{K}(0,\nu)=\sum\limits_{k=0
}^{K}\mathbf{E}[\gamma^{2}\|\nu(k)\|^{2}-\|z(k)\|^{2}]\geq 0$, which implies
$\|L_{K}\|\leq\gamma$. Following the line of Lemma 3 of \cite{ref11},
we can further show $\|L_{K}\|<\gamma$ with the detail omitted.

Similarly to above process, we have from Lemma \ref{Le:2.3} that
\begin{eqnarray*}
&&\min_{\nu\in
l^{2}_{\omega}(N_{K},R^{l})}J^{K}(x_{k_{0}},\nu)=J^{K}(x_{k_{0}},\nu^{*})\nonumber\\
&=&[\mathbf{E}x_{k_{0}}]^{T}Q_{1}(k_{0})[\mathbf{E}x_{k_{0}}] \leq  J^{K}(x_{k_{0}},0)\nonumber\\
&=&-\sum\limits_{k=k_{0} }^{K}\mathbf{E}[\|z(k)\|^{2}] \leq  0
\end{eqnarray*}
for arbitrary $x_{k_{0}}\in R^{n}$. So $Q_{1}(k)\leq 0$, $k\in
N_{K}$. Theorem \ref{th:2.4} is proved.\rulex

\begin{remark}\label{re:2.5}
Theorem \ref{th:2.4} is only a sufficient but not a
necessary condition for  $\|L_{K}\|<\gamma$, which is different from
 classical discrete-time stochastic systems \cite{ref11}. For the constrained
backward difference equation  (\ref{E:2.6}), due to $P(K+1)=Q(K+1)=0$,
$H(P(K+1))>0$ and $\tilde{H}(P(K+1),Q(K+1))>0$, we can get a  unique
solution  $(P(K), Q(K))$. Similarly, we can compute $(P(K-1),
Q(K-1))$ if $H(P(K))>0$ and $\tilde{H}(P(K),Q(K))>0$. The equation
(\ref{E:2.6}) can be solved backwardly for ever if and only if
$H(P(t+1))>0, \tilde{H}(P(t+1),Q(t+1))>0$ for $t=k-2, k-3, ..., 0$.
However, $\|L_{K}\|<\gamma$ does not necessarily imply  $H(P(t+1))>0$ and
$\tilde{H}(P(t+1),Q(t+1))>0$ simultaneously, so Theorem \ref{th:2.4}
is only a sufficient condition, the solvability of
(\ref{E:2.6})  merits further study.
\end{remark}

\section{Mean-field Stochastic LQ Control}

Consider the following discrete-time stochastic difference equation
\begin{equation}
\label{E:3.1}
\begin{cases}
 \  \bar{x}(k+1)= A_{1}(k)\bar{x}(k)+\tilde{A}_{1}(k)\mathbf{E}\bar{x}(k)+F_{1}(k)\bar{u}(k)+[B_{1}(k)\bar{x}(k)+\tilde{B}_{1}(k)\mathbf{E}\bar{x}(k)]\omega(k),\\
 \ \bar{z}(k)=\left[\begin{array}{c}
{\it \Phi}_{1}(k)\bar{x}(k) \\
{\it \Psi}_{1}(k)\bar{u}(k)
\end{array}\right] ,& \\
{\it \Psi}_{1}^{T}(k){\it \Psi}_{1}(k)=I,\q \bar{x}(0)=\bar{x}_{0}, k\in N_{K},
\end{cases}
\end{equation}
where $\bar{u}(k)\in l^{2}_{\omega}(N_{K},R^{q})$ is the control input. The associated cost function is
\[
J^{K}(\bar{x}_{0},\bar{u})=\d\sum\limits_{k=0
}^{K}\mathbf{E}[\|\bar{z}(k)\|^{2}].
\]
Similarly to the proof of Theorem \ref{th:2.4}, it is easy to obtain

\begin{theorem}\label{th:3.1} (LQ control) For the  mean-field type
stochastic system (\ref{E:3.1}), there exists $\bar{u}^{*}\in
l^{2}_{\omega}(N_{K},R^{q})$ such that $\mathop {\min
}\limits_{\bar{u}\in
l^{2}_{\omega}(N_{K},R^{q})}J^{K}(\bar{x}_{0},\bar{u})=J^{K}(\bar{x}_{0},\bar{u}^{*})
=[\mathbf{E}\bar{x}_{0}]^{T}\tilde{Q}_{1}(0)[\mathbf{E}\bar{x}_{0}]\geq 0$ and
$\tilde{Q}_{1}(k)\geq 0$ if the following backward difference
equation
\begin{equation}
\label{E:3.2}
\begin{cases}
 \  \tilde{P}_{1}(k)=L_{1}(\tilde{P}_{1}(k+1))-G_{1}(\tilde{P}_{1}(k+1))H_{1}(\tilde{P}_{1}(k+1))^{-1}G_{1}(\tilde{P}_{1}(k+1))^{T}, &\\
 \  \tilde{Q}_{1}(k)=\tilde{L}_{1}(\tilde{P}_{1}(k+1),\tilde{Q}_{1}(k+1))-\tilde{G}_{1}(\tilde{P}_{1}(k+1),\tilde{Q}_{1}(k+1))&\\
 \ \hspace{1.4cm}\times\tilde{H}_{1}(\tilde{P}_{1}(k+1),\tilde{Q}_{1}(k+1))^{-1}\tilde{G}_{1}(\tilde{P}_{1}(k+1),\tilde{Q}_{1}(k+1))^{T} , & \\
 \ \tilde{P}_{1}(K+1)=\tilde{Q}_{1}(K+1)=0, &\\
 \ H_{1}(\tilde{P}_{1}(k+1))>0, \q \tilde{H}_{1}(\tilde{P}_{1}(k+1), \tilde{Q}_{1}(k+1))>0
  \end{cases}
\end{equation}
has a unique solution $(\tilde{P}_{1}(k),\tilde{Q}_{1}(k))$ with
$k\in N_{K}$,
where
\begin{eqnarray*}
&&\bar{u}^{*}(k)=-H_{1}(\tilde{P}_{1}(k+1))^{-1}G_{1}(\tilde{P}_{1}(k+1))^{T}\bar{x}(k)+[H_{1}(\tilde{P}_{1}(k+1))^{-1}G_{1}
(\tilde{P}_{1}(k+1))^{T}
\nonumber\\
&&\hspace{1.4cm}-\tilde{H}_{1}(\tilde{P}_{1}(k+1),\tilde{Q}_{1}(k+1))^{-1}\tilde{G}_{1}(\tilde{P}_{1}(k+1),\tilde{Q}_{1}(k+1))^{T}]E\bar{x}(k),
\nonumber\\
&&L_{1}(\tilde{P}_{1}(k+1))=A_{1}(k)^{T}\tilde{P}_{1}(k+1)A_{1}(k)
+B_{1}(k)^{T}\tilde{P}_{1}(k+1)B_{1}(k)+{\it \Phi}_{1}(k)^{T}{\it \Phi}_{1}(k),\nonumber\\
&&G_{1}(\tilde{P}_{1}(k+1))=A_{1}(k)^{T}\tilde{P}_{1}(k+1)F_{1}(k),
 \q H_{1}(\tilde{P}_{1}(k+1))=I_{q}+F_{1}(k)^{T}\tilde{P}_{1}(k+1)F_{1}(k),\nonumber\\
&&\mathbb{A}_{1}(k)=A_{1}(k)+\tilde{A}_{1}(k),\q \mathbb{B}_{1}(k)=B_{1}(k)+\tilde{B}_{1}(k),\nonumber\\
&&\tilde{L}_{1}(\tilde{P}_{1}(k+1),\tilde{Q}_{1}(k+1))=\mathbb{A}_{1}(k)^{T}\tilde{Q}_{1}(k+1)\mathbb{A}_{1}(k)
+\mathbb{B}_{1}(k)^{T}\tilde{P}_{1}(k+1)\mathbb{B}_{1}(k)+{\it \Phi}_{1}(k)^{T}{\it \Phi}_{1}(k),
\nonumber\\
&&\tilde{G}_{1}(\tilde{P}_{1}(k+1),\tilde{Q}_{1}(k+1))
=\mathbb{A}_{1}(k)^{T}\tilde{Q}_{1}(k+1)F_{1}(k),
\nonumber\\
&&\tilde{H}_{1}(\tilde{P}_{1}(k+1),\tilde{Q}_{1}(k+1))=I_{q}+F_{1}(k)^{T}\tilde{Q}_{1}(k+1)F_{1}(k).
\end{eqnarray*}
\end{theorem}

\section{Main Results}
We first define the finite-time $H_2/H_\infty$ control as follows:

\begin{definition}\label{df:4.1} Consider the  controlled
stochastic system (\ref{E:1.1}) with $k\in N_{K}$,
where $u(k)\in l^{2}_{\omega}(N_{K},R^{q})$ is the control input.
Given $0<K<\infty$ and the disturbance attenuation level $\gamma>0$,
if existing, a state feedback control
$u^{*}(k)=U(k)x(k)+\tilde{U}(k)\mathbf{E}x(k)=U(k)(x(k)-\mathbf{E}x(k))+(U(k)+\tilde{U}(k))\mathbf{E}x(k)\in
l^{2}_{\omega}(N_{K},R^{q})$, such that

1) For the closed-loop system
\begin{equation}
\label{E:4.2}
\begin{cases}
 \  x(k+1)= (A(k)+F_{1}(k)U(k))x(k)+(\tilde{A}(k)+F_{1}(k)\tilde{U}(k))\mathbf{E}x(k)&\\
\ \hspace{1.6cm}+B(k)\nu(k)+\tilde{B}(k)\mathbf{E}\nu(k)&\\
\ \hspace{1.6cm}+[C(k)x(k)+\tilde{C}(k)\mathbf{E}x(k)+D(k)\nu(k)+\tilde{D}(k)\mathbf{E}\nu(k)]\omega(k), & \\
 \ z(k)=\left[\begin{array}{c}
{\it \Phi}(k)x(k) \\
{\it \Psi}(k)(U(k)x(k)+\tilde{U}(k)\mathbf{E}x(k))
\end{array}\right], & \\
{\it \Psi}^{T}(k){\it \Psi}(k)=I, \q  x(0)=x_{0}\in R^{n}, &
\end{cases}
\end{equation}
the following
\[\begin{split}
\|L_{K}\|=\mathop {\sup }\limits_{\nu\in
l^{2}_{\omega}(N_{K},R^{l}),\nu\neq 0,x_{0}=0}
\frac{\|z(k)\|_{l^{2}_{\omega}(N_{K},R^{m})}}{\|\nu(k)\|_{l^{2}_{\omega}(N_{K},R^{l})}}<\gamma
\end{split}\]
holds.

2) When the worst case disturbance
$\nu^{*}(k)=V(k)x(k)+\tilde{V}(k)\mathbf{E}x(k)$,  if existing, is
implemented in (\ref{E:1.1}), $u^{*}(k)$ minimizes the output energy
$
J^{K}_{2}(u,\nu^{*})=\|z(k)\|^2_{{l^{2}_{\omega}(N_{K},R^{m})}}.
$
\end{definition}

If $(u^*,\nu^{*})$ exists, we also say that the finite horizon
$H_2/H_\infty$ control of mean-field type is solvable. Before
presenting the main result, we introduce four coupled matrix-valued
equations  as
\begin{equation}
\label{E:4.3}
\begin{cases}
 \  P_{1}(k)=(A(k)+F_{1}(k)U(k))^{T}P_{1}(k+1)(A(k)+F_{1}(k)U(k))\\
 \ \hspace{1.3cm}+C(k)^{T}P_{1}(k+1)C(k)-{\it \Phi}(k)^{T}{\it \Phi}(k)-U(k)^{T}U(k)\\
 \ \hspace{1.3cm}-G_{u}(P_{1}(k+1))H(P_{1}(k+1))^{-1}G_{u}(P_{1}(k+1))^{T}, &\\
 \  Q_{1}(k)=[\mathbb{A}(k)+F_{1}(k)\mathbb{U}(k)]^{T}Q_{1}(k+1))[\mathbb{A}(k)+F_{1}(k)\mathbb{U}(k)]\\
 \  \hspace{1.3cm}+\mathbb{C}(k)^{T}P_{1}(k+1)\mathbb{C}(k)-{\it \Phi}(k)^{T}{\it \Phi}(k)-\mathbb{U}(k)^{T}\mathbb{U}(k)\\
  \  \hspace{1.3cm}-\tilde{G}_{u}(P_{1}(k+1),Q_{1}(k+1))
  \tilde{H}(P_{1}(k+1),Q_{1}(k+1))^{-1}(\cdot\cdot\cdot)^{T}, & \\
 \ P_{1}(K+1)=Q_{1}(K+1)=0, &\\
 \ H(P_{1}(k+1))>0, \q \tilde{H}(P_{1}(k+1), Q_{1}(k+1))>0, &
  \end{cases}
\end{equation}
\begin{equation}
\label{E:4.4}
\begin{cases}
 V(k)=\lefteqn{-H(P_{1}(k+1))^{-1}G_{u}(P_{1}(k+1))^{T},}\\
\mathbb{V}(k)=V(k)+\tilde{V}(k)=\tilde{H}(P_{1}(k+1),Q_{1}(k+1))^{-1}\tilde{G}_{u}(P_{1}(k+1),Q_{1}(k+1))^{T},
 \end{cases}
\end{equation}
\begin{equation}
\label{E:4.5}
\begin{cases}
 \  \tilde{P}_{1}(k)=[A(k)+B(k)V(k)]^{T}\tilde{P}_{1}(k+1)[A(k)+B(k)V(k)]\\
 \ \hspace{1.3cm}+[C(k)+D(k)V(k)]^{T}\tilde{P}_{1}(k+1)[C(k)+D(k)V(k)]+{\it \Phi}(k)^{T}{\it \Phi}(k)+I_{q} \\
\ \hspace{1.3cm}-G_{\nu}(\tilde{P}_{1}(k+1))H_{1}(\tilde{P}_{1}(k+1))^{-1}G_{\nu}(\tilde{P}_{1}(k+1))^{T}, &\\
  \  \tilde{Q}_{1}(k)=[\mathbb{A}(k)+\mathbb{B}(k)\mathbb{V}(k)]^{T}\tilde{Q}_{1}(k+1))
  [\mathbb{A}(k)+\mathbb{B}(k)\mathbb{V}(k)]\\
 \ \hspace{1.3cm}+[\mathbb{C}(k)+\mathbb{D}(k)\mathbb{V}(k)]^{T}\tilde{P}_{1}(k+1)
 [\mathbb{C}(k)+\mathbb{D}(k)\mathbb{V}(k)]+{\it \Phi}(k)^{T}{\it \Phi}(k)+I_{q},\\
\hspace{1.3cm}-\tilde{G}_{\nu}(\tilde{P}_{1}(k+1),\tilde{Q}_{1}(k+1))
\tilde{H}_{1}(\tilde{P}_{1}(k+1),\tilde{Q}_{1}(k+1))^{-1}
(\cdot\cdot\cdot)^{T}, &\\
  \ \tilde{P}_{1}(K+1)=\tilde{Q}_{1}(K+1)=0, & \\
  \ H_{1}(\tilde{P}_{1}(k+1))>0, \q \tilde{H}_{1}(\tilde{P}_{1}(k+1), \  \tilde{Q}_{1}(k+1))>0, &
  \end{cases}
\end{equation}
\begin{equation}
\label{E:4.6}
\begin{cases}
 U(k)=\lefteqn{-H_{1}(\tilde{P}_{1}(k+1))^{-1}G_{\nu}(\tilde{P}_{1}(k+1))^{T},}\\
 \mathbb{U}(k)=U(k)+\tilde{U}(k)=\tilde{H}_{1}(\tilde{P}_{1}(k+1),\tilde{Q}_{1}(k+1))^{-1}
\tilde{G}_{\nu}(\tilde{P}_{1}(k+1),\tilde{Q}_{1}(k+1))^{T},
\end{cases}
\end{equation}
where
\ay
\begin{eqnarray*}
&&G_{u}(P_{1}(k+1))=(A(k)+F_{1}(k)U(k))^{T}P_{1}(k+1)B(k)+C(k)P_{1}(k+1)D(k),
\nonumber\\
&&G_{\nu}(\tilde{P}_{1}(k+1))=[A(k)+B(k)V(k)]^{T}\tilde{P}_{1}(k+1)F_{1}(k),
\nonumber\\
&&\tilde{G}_{u}(P_{1}(k+1),Q_{1}(k+1))=[\mathbb{A}(k)+F_{1}(k)\mathbb{U}(k)]^{T}Q_{1}(k+1)\mathbb{B}(k)
+\mathbb{C}(k)^{T}P_{1}(k+1)\mathbb{D}(k),
\nonumber\\
&&\tilde{G}_{\nu}(\tilde{P}_{1}(k+1),\tilde{Q}_{1}(k+1))=
[\mathbb{A}(k)+\mathbb{B}(k)\mathbb{V}(k)]^{T}\tilde{Q}_{1}(k+1)F_{1}(k).
\end{eqnarray*}

Our main result in this section is given by the following theorem:

\begin{theorem}\label{th:4.2}For  a  given disturbance attenuation lever
$\gamma>0$, the finite horizon $H_2/H_\infty$ control of mean-field
type system has  the  solution $(u^{*}(k), \nu^{*}(k))$ as
$$u^{*}(k)=U(k)x(k)+\tilde{U}(k)\mathbf{E}x(k),\hspace{20pt}\nu^{*}(k)=V(k)x(k)+\tilde{V}(k)\mathbf{E}x(k) $$
with $U(k), \tilde{U}(k)\in R^{q\times n}$ and $V(k),
\tilde{V}(k)\in R^{l\times n}$ being matrix-valued functions and
$Q_{1}(k)\leq 0,\tilde{Q}_{1}(k)\geq 0$, if the  matrix-valued
equations (\ref{E:4.3})-(\ref{E:4.6}) have the  solution
$(P_{1}(k)$,$Q_{1}(k)$; $\tilde{P}_{1}(k)$, $\tilde{Q}_{1}(k)$;
$U(k)$, $\tilde{U}(k)$; $V(k)$, $\tilde{V}(k))$ with $k\in N_{K}$.
\end{theorem}

\proof With the solution
$(P_{1}(k),Q_{1}(k);\tilde{P}_{1}(k),\tilde{Q}_{1}(k);U(k),
\tilde{U}(k);V(k), \tilde{V}(k))$ to  the  equations
(\ref{E:4.3})-(\ref{E:4.6}), we can construct
$u^{*}(k)=U(k)x(k)+\tilde{U}(k)\mathbf{E}x(k)$ and substitute $u^{*}(k)$ into
system (\ref{E:1.1}), then system (\ref{E:4.2}) is obtained. By Theorem
\ref{th:2.4} and (\ref{E:4.3}), it yields that $\|L_{K}\|<\gamma$. Keeping
(\ref{E:4.3}) in mind, by the technique of completing squares and
Lemma 2.2, we immediately get $Q_{1}(k)\leq 0$ and
\begin{eqnarray*}
J^{K}_{1}(u^{*},\nu)&=&\d\sum\limits_{k=0
}^{K}\mathbf{E}[\gamma^{2}\|\nu(k)\|^{2}-\|z(k)\|^{2}]\nonumber\\
&=&\d\sum\limits_{k=0 }^{K} \mathbf{E} \left\{\left[\begin{array}{c}
x(k)-\mathbf{E}x(k) \\
\nu(k)-\mathbf{E}\nu(k)
\end{array}\right]^{T}M_{0}(P_{1}(k))
\left[\begin{array}{c}
x(k)-\mathbf{E}x(k) \\
\nu(k)-\mathbf{E}\nu(k)
\end{array}\right]\right\}\nonumber\\
&&+\d\sum\limits_{k=0 }^{K} \left[\begin{array}{c}
\mathbf{E}x(k) \\
\mathbf{E}\nu(k)
\end{array}\right]^{T}S_{0}(P_{1}(k),Q_{1}(k))
\left[\begin{array}{c}
\mathbf{E}x(k) \\
\mathbf{E}\nu(k)
\end{array}\right]
+[\mathbf{E}x_{0}]^{T}Q_{1}(0)[\mathbf{E}x_{0}]
\nonumber\\
&=&\d\sum\limits_{k=0 }^{K}
\mathbf{E}\{[(\nu(k)-\mathbf{E}\nu(k))-(\nu^{*}(k)-\mathbf{E}\nu^{*}(k))]^{T}
H(P_{1}(k+1))[\cdots]\}
\nonumber\\
&&+[\mathbf{E}x_{0}]^{T}Q_{1}(0)[\mathbf{E}x_{0}]+\sum\limits_{k=0 }^{K}
[\mathbf{E}\nu(k)-\mathbf{E}\nu^{*}(k)]^{T}\tilde{H}(P_{1}(k+1),Q_{1}(k+1))[\cdots]
\nonumber\\
&\geq&J^{K}_{1}(u^{*},\nu^{*})=[\mathbf{E}x_{0}]^{T}Q_{1}(0)[\mathbf{E}x_{0}].
\end{eqnarray*}
So, we see that $\nu^{*}(k)=V(k)x(k)+\tilde{V}(k)\mathbf{E}x(k)$ with
$(V(k),\tilde{V}(k))$ given by (\ref{E:4.4}) is the worse case
disturbance, where
\ay
\begin{eqnarray*}
&&M_{0}(P_{1}(k))=
\left[\begin{array}{cc}
L_{0}(P_{1}(k+1)) &  G_{u}(P_{1}(k+1))     \\
G_{u}(P_{1}(k+1))^{T} &  H(P_{1}(k+1))
 \end{array}\right],\nonumber\\
&&S_{0}(P_{1}(k),Q_{1}(k))=
\left[\begin{array}{cc}
\tilde{L}_{0}(\tilde{P}_{1}(k+1),\tilde{Q}_{1}(k+1)) &  \tilde{G}_{u}(\tilde{P}_{1}(k+1),\tilde{Q}_{1}(k+1))     \\
 \tilde{G}_{u}(\tilde{P}_{1}(k+1),\tilde{Q}_{1}(k+1))^{T} &  \tilde{H}(\tilde{P}_{1}(k+1),\tilde{Q}_{1}(k+1))
 \end{array}\right]
\end{eqnarray*}
with
\ay
\begin{eqnarray*}
&&L_{0}(P_{1}(k+1))=-P_{1}(k)+(A(k)+F_{1}(k)U(k))^{T}P_{1}(k+1))(A(k)+F_{1}(k)U(k))
\nonumber\\
  &&\hspace{2.5cm}+C(k)^{T}P_{1}(k+1)C(k)-{\it \Phi}(k)^{T}{\it \Phi}(k)-U(k)^{T}U(k),
\nonumber\\
&&\tilde{L}_{0}(\tilde{P}_{1}(k+1),\tilde{Q}_{1}(k+1))=-\tilde{Q}_{1}(k)+[\mathbb{A}(k)+F_{1}(k)\mathbb{U}(k)]^{T}\tilde{Q}_{1}(k+1))[\cdots]
\nonumber\\
&&\hspace{4.4cm}+\mathbb{C}(k)^{T}\tilde{P}_{1}(k+1)\mathbb{C}(k)-{\it \Phi}(k)^{T}{\it \Phi}(k)-\mathbb{U}(k)^{T}\mathbb{U}(k).
\end{eqnarray*}
Similarly,  Theorem \ref{th:3.1} and (\ref{E:4.5})  yield
$\tilde{Q}_{1}(k)\geq 0$ and
\[\begin{split}
J^{K}_{2}(u,\nu^{*})=
\lefteqn{\d\sum\limits_{k=0 }^{K}\mathbf{E}[\|z(k)\|^{2}]
}
\nonumber\\
=&[\mathbf{E}x_{0}]^{T}\tilde{Q}_{1}(0)[\mathbf{E}x_{0}]+\d\sum\limits_{k=0 }^{K}
[\mathbf{E}u(k)-\mathbf{E}u^{*}(k)]^{T}\tilde{H}_{1}(\tilde{P}_{1}(k+1),\tilde{Q}_{1}(k+1))[\cdots]
\nonumber\\
&+\d\sum\limits_{k=0 }^{K}
\mathbf{E}\{[(u(k)-\mathbf{E}u(k))-(u^{*}(k)-\mathbf{E}u^{*}(k))]^{T}H_{1}(\tilde{P}_{1}(k+1))[\cdots]\}
\nonumber\\
\geq&J^{K}_{2}(u^{*},\nu^{*})=[\mathbf{E}x_{0}]^{T}\tilde{Q}_{1}(0)[\mathbf{E}x_{0}].
\end{split}\]
Therefore, $(u^{*},\nu^{*})$ solve the mean-field $H_2/H_\infty$ control problem of system (\ref{E:1.1}), and the proof is complete.\rulex

\begin{remark}\label{re:4.3}   For the matrix-valued equations
(\ref{E:4.3})-(\ref{E:4.6}),  from  $P_{1}(K+1)=Q_{1}(K+1)=0,
\tilde{P}_{1}(K+1)=\tilde{Q}_{1}(K+1)=0$, we know $H(P_{1}(K+1))>0,
\tilde{H}(P_{1}(K+1),Q_{1}(K+1))>0,
H_{1}(\tilde{P}_{1}(K+1))>0,\tilde{H}_{1}(\tilde{P}_{1}(K+1),\tilde{Q}_{1}(K+1))>0$.
Accordingly,   $(U(K),V(K))$ and $(\tilde{U}(K),\tilde{V}(K))$ can
be computed by  the matrix equations (\ref{E:4.4}) and (\ref{E:4.6}),
then $(P_{1}(K), Q_{1}(K)\leq 0)$ and $(\tilde{P}_{1}(K),
\tilde{Q}_{1}(K)\geq 0)$ can be obtained by (\ref{E:4.3}) and
(\ref{E:4.5}). The backward  recursion can proceed   if and only if
$H(P_{1}(t))>0, \tilde{H}(P_{1}(t),Q_{1}(t))>0,
H_{1}(\tilde{P}_{1}(t))>0,\tilde{H}_{1}(\tilde{P}_{1}(t),\tilde{Q}_{1}(t)>0$
for $t=k-1,k-2, k-3, ..., 0$. In   (\ref{E:1.1}), if $\tilde A\equiv
0$,   $\tilde B\equiv 0$,  $\tilde C\equiv 0$,  $\tilde D\equiv 0$,
the solvability of the finite horizon $H_2/H_\infty$ control is
equivalent to that of  the matrix-valued equations
(\ref{E:4.3})-(\ref{E:4.6}); see \cite{ref11}.  However, for
(\ref{E:1.1}),  the solvability condition of (\ref{E:4.3})-(\ref{E:4.6})
remains unsolved at present stage.
\end{remark}

\begin{remark}\label{re:4.4}
In this paper, the disturbance attenuation level $\gamma>0$ is given
in advance, the definition of our  mixed $H_2/H_\infty$ control
arises from the classical work [2]. If $\gamma>0$ is not
predetermined or in other words, we have to   select $\gamma>0$  to
ensure a good trade off  between the two contradictory objectives
$H_2$ optimization and the $H_\infty$ optimal disturbance level,
this is another issue called multi-objective  $H_2/H_\infty$
control; see \cite{kzhou,linchen}.
\end{remark}

\begin{remark}\label{re:4.4}  Mean-field
stochastic systems have been used to mean-variance portfolio
selection \cite{Hafayed}, Social optima \cite{ref27} and large
population systems \cite{ref26}, where in these works, the external
disturbance is not considered in mathematical modeling. Generally
speaking, in a real world, the exogenous influence always exists.
For example, in a financial market, the stock price is subject to
unexpected disaster and  political strategy, which can be
represented by $\omega(\cdot)$.  So, it is expected that what we
have obtained may be useful in mathematical finance and other
practical fields, which motivates us to do this research.
\end{remark}

\section{Algorithm and Numerical Example}

If the matrix-valued equations  (\ref{E:4.3})-(\ref{E:4.6})  are
solvable, they can  be solved recursively as follows:

\begin{itemize}

\item [i)] Let $k=K$, then
$H(P_{1}(K+1)),\tilde{H}(P_{1}(K+1),Q_{1}(K+1)),
H_{1}(\tilde{P}_{1}(K+1)),\tilde{H}_{1}(\tilde{P}_{1}(K+1),\tilde{Q}_{1}(K+1)),
 $ can be computed by $P_{1}(K+1)=Q_{1}(K+1)=0, \tilde{P}_{1}(K+1)=\tilde{Q}_{1}(K+1)=0$.

\item [ii)]   If $H(P_{1}(K+1))>0,\tilde{H}(P_{1}(K+1),Q_{1}(K+1))>0,
H_{1}(\tilde{P}_{1}(K+1))>0,\tilde{H}_{1}(\tilde{P}_{1}(K+1),\tilde{Q}_{1}(K+1))>0,
 $, calculate $H(P_{1}(K+1))^{-1},  \tilde{H}(P_{1}(K+1),Q_{1}(K+1))^{-1},
H_{1}(\tilde{P}_{1}(K+1))^{-1}, \tilde{H}_{1}(\tilde{P}_{1}(K+1), \tilde{Q}_{1}(K+1))^{-1}$.

\item [iii)] Solving the matrix equations (\ref{E:4.4}) and (\ref{E:4.6})
to obtain $(U(K),\tilde{U}(K),V(K),\tilde{V}(K))$.

\item [iv)] Substitute the obtained $(U(K),\tilde{U}(K))$ into the  matrix
equation (\ref{E:4.3}) and $(V(K),\tilde{V}(K))$ into (\ref{E:4.5}),
then $(P_{1}(K), Q_{1}(K)\leq 0, \tilde{P}_{1}(K),
\tilde{Q}_{1}(K)\geq 0)$ are available.

\item [v)] Repeat the above procedures,
$(U(k),\tilde{U}(k),V(k),\tilde{V}(k))$ and $(P_{1}(k), Q_{1}(k),
\tilde{P}_{1}(k), \tilde{Q}_{1}(k))$ can be computed recursively for
$k=K-1,K-2,K-3,\cdots,0$.

\end{itemize}

Next, we present a two-step  numerical example to show the detail
and efficiency of the above algorithm.

\begin{example} In  system (\ref{E:1.1}),  set  $K=2, \gamma=0.8$.
the parameters of system (\ref{E:1.1}) is given   in Table 1.
According to the above algorithm, we can check the existence of the
solutions of the coupled matrix-valued equations
(\ref{E:4.3})-(\ref{E:4.6}) and obtain them backward. Table 2
illustrates the solutions.
\begin{center}
{\small
    \begin{center}
    \centerline{\small {\bf Table 1}~~Parameters of system (\ref{E:1.1})}\vskip 1mm
    \label{Tab:MSE}

    {\small
    \begin{tabular*}{11.5cm}{ccccc}
        \toprule
          &time &$k=2$ & $k=1$ &$k=0$   \\[0.8ex]
        \midrule
        &$A(k)$
        &$\left[ \begin{array}{*{2}{p{1.1cm}}} 0.1500 & 0.1000 \\0.2000 & 0.1500 \end{array}\right]$
        &$\left[ \begin{array}{cc} 0.1000 & 0.0800 \\0.1800 & 0.1200 \end{array}\right]$
        &$\left[ \begin{array}{*{2}{p{1.1cm}}} 0.0500 & 0.1500 \\0.2500 &0.3500 \end{array}\right]$  \\[3ex]
        &$\tilde{A}(k)$
        &$\left[ \begin{array}{*{2}{p{1.1cm}}} 0.1500 & 0.1500 \\0.2500 & 0.1000 \end{array}\right]$
        &$\left[ \begin{array}{cc} 0.1000 & 0.1200 \\0.2200 & 0.0800 \end{array}\right]$
        &$\left[ \begin{array}{*{2}{p{1.1cm}}} 0.0500 & 0.2500 \\0.3500 & 0.2000 \end{array}\right]$  \\[3ex]
        &$B(k)$
        &$\left[ \begin{array}{*{2}{p{1.1cm}}} 0.1500 & 0.2000 \\0.2000 & 0.3000 \end{array}\right]$
        &$\left[ \begin{array}{cc} 0.1000 & 0.1800 \\0.2000 & 0.2800 \end{array}\right]$
        &$\left[ \begin{array}{*{2}{p{1.1cm}}} 0.0500 & 0.1000 \\0.1000 & 0.2000 \end{array}\right]$   \\[3ex]
        &$\tilde{B}(k)$
        &$\left[ \begin{array}{*{2}{p{1.1cm}}} 0.1500 & 0.2500 \\0.3000 & 0.1000 \end{array}\right]$
        &$\left[ \begin{array}{cc} 0.1000 & 0.2000 \\0.2500 & 0.0800 \end{array}\right]$
        &$\left[ \begin{array}{*{2}{p{1.1cm}}} 0.0500 & 0.1500 \\0.2000 & 0.2000 \end{array}\right]$  \\[3ex]
        &$C(k)$
        &$\left[ \begin{array}{*{2}{p{1.1cm}}} 0.1500 & 0.1500 \\0.2000 & 0.1500 \end{array}\right]$
        &$\left[ \begin{array}{cc} 0.1000 & 0.1200 \\0.1800 & 0.1000 \end{array}\right]$
        &$\left[ \begin{array}{*{2}{p{1.1cm}}} 0.0500 & 0.2500 \\0.1000 & 0.2500 \end{array}\right]$  \\[3ex]
        &$\tilde{C}(k) $
        &$\left[ \begin{array}{*{2}{p{1.1cm}}} 0.1500 & 0.1000 \\0.2000 & 0.1000 \end{array}\right]$
        &$\left[ \begin{array}{cc} 0.1000 & 0.0800 \\0.1500 & 0.0800 \end{array}\right]$
        &$\left[ \begin{array}{*{2}{p{1.1cm}}} 0.0500 & 0.1500 \\0.2500 & 0.1500 \end{array}\right]$  \\[3ex]
        &$D(k)$
        &$\left[ \begin{array}{*{2}{p{1.1cm}}}0.1500 & 0.1000 \\0.1500 & 0.2000 \end{array}\right]$
        &$\left[ \begin{array}{cc} 0.1000 & 0.1200 \\0.2000 & 0.1800 \end{array}\right]$
        &$\left[ \begin{array}{*{2}{p{1.1cm}}} 0.0500 & 0.1800 \\0.1500 & 0.2800 \end{array}\right]$  \\[3ex]
        &$\tilde{D}(k)$
        &$\left[ \begin{array}{*{2}{p{1.1cm}}} 0.1500 & 0.1500 \\0.2000 & 0.2500 \end{array}\right]$
        &$\left[ \begin{array}{cc} 0.1000 & 0.4000 \\0.1000 & 0.1500 \end{array}\right]$
        &$\left[ \begin{array}{*{2}{p{1.1cm}}} 0.0500 & 0.3000 \\0.2500 & 0.3500 \end{array}\right]$  \\[3ex]
        &$F_{1}(k)$
        &$\left[ \begin{array}{*{2}{p{1.1cm}}} 0.1500 & 0.2000 \\0.1500 & 0.2000 \end{array}\right]$
        &$\left[ \begin{array}{cc} 0.1000 & 0.3000 \\0.2500 & 0.1000 \end{array}\right]$
        &$\left[ \begin{array}{*{2}{p{1.1cm}}} 0.0500 & 0.2500 \\0.3500 & 0.3000 \end{array}\right]$  \\[3ex]
        &$\Phi(k)$
        &$\left[ \begin{array}{*{2}{p{1.1cm}}} 0.1500 & 0.1500 \\0.2000 & 0.3000 \end{array}\right]$
        &$\left[ \begin{array}{cc} 0.1000 & 0.2500 \\0.1000 & 0.2000 \end{array}\right]$
        &$\left[ \begin{array}{*{2}{p{1.1cm}}} 0.0500 & 0.1000 \\0.3000 & 0.2000 \end{array}\right]$  \\[3ex]
        &$\Psi(k)$
        &$\left[ \begin{array}{*{2}{p{1.1cm}}}  0.6000 & -0.8000 \\0.8000 & 0.6000 \end{array}\right]$
        &$\left[ \begin{array}{cc} 1.0000 & 0.0000 \\0.0000 & 1.0000 \end{array}\right]$
        &$\left[ \begin{array}{*{2}{p{1.1cm}}} 0.8000 & 0.6000 \\0.6000 & -0.8000 \end{array}\right]$  \\[3ex]
        \bottomrule
    \end{tabular*}}
\end{center}}
\end{center}

\begin{center}
{\small
    \begin{center}
    \centerline{\small {\bf Table 2}~~Solutions for (\ref{E:4.3})-(\ref{E:4.6})}\vskip 1mm
    \label{Tab:MSE1}
    {\small
    \begin{tabular*}{14.5cm}{ccccc}
        \toprule
          &time &$k=2$ & $k=1$ &$k=0$   \\[0.8ex]
        \midrule
        &$H(P_{1}(k+1))$
        &$\left[ \begin{array}{*{2}{p{1.15cm}}} \hspace{0.1cm}0.6400 & \hspace{0.6cm}0 \\ \hspace{0.6cm}0 & \hspace{0.1cm}0.6400  \end{array}\right]$
        &$\left[ \begin{array}{cc} \ \ 0.6232 & -0.0210 \\-0.0210 & \ \ \ 0.6127 \end{array}\right]$
        &$\left[ \begin{array}{cc} \ \ 0.6346 & -0.0112 \\-0.0112 & \ \ \ 0.6166 \end{array}\right]$  \\[3ex]
        &$\begin{array}{cc} $$\tilde{H}(P_{1}(k+1),Q_{1}(k+1))$$  \end{array}$
        &$\left[ \begin{array}{*{2}{p{1.15cm}}} \hspace{0.1cm}0.6400 & \hspace{0.6cm}0 \\ \hspace{0.6cm}0 & \hspace{0.1cm}0.6400  \end{array}\right]$
        &$\left[ \begin{array}{cc} \ \ 0.5773 & -0.0790 \\-0.0790 & \ \ 0.5364 \end{array}\right]$
        &$\left[ \begin{array}{cc} \ \ 0.5950 & -0.0801 \\-0.0801 & \ \ \ 0.4948 \end{array}\right]$   \\[3ex]
        &$H_{1}(\tilde{P}_{1}(k+1))$
        &$\left[ \begin{array} {*{2}{p{1.15cm}}} \hspace{0.6cm}1 & \hspace{0.6cm}0 \\ \hspace{0.6cm}0 & \hspace{0.6cm}1 \end{array}\right]$
        &$\left[ \begin{array}{cc} \ \  1.0843  & \ \ \ 0.0667\\ \ \ 0.0667 &  \ \ \ 1.1117 \end{array}\right]$
        &$\left[ \begin{array}{cc}  \ \ 1.1489 & \ \ \ 0.1480 \\ \ \ 0.1480 & \ \ \ 1.1925 \end{array}\right]$  \\[3ex]
        &$\begin{array}{cc} $$\tilde{H}(\tilde{P}_{1}(k+1), \tilde{Q}_{1}(k+1))$$  \end{array}$
        &$\left[ \begin{array}{*{2}{p{1.15cm}}} \hspace{0.6cm}1 & \hspace{0.6cm}0 \\ \hspace{0.6cm}0 & \hspace{0.6cm}1 \end{array}\right]$
        &$\left[ \begin{array}{cc} \ \ 1.0843 & \ \ \ 0.0667 \\  \ \ 0.0667 &   \ \ 1.1117 \end{array}\right]$
        &$\left[ \begin{array}{cc} \ \ 1.1902 & \ \ \ 0.2139 \\ \ \ 0.2139 & \ \ \ 1.2985 \end{array}\right]$  \\[3ex]
        &$U(k)$
        &$\left[ \begin{array}{*{2}{p{1.15cm}}} \hspace{0.6cm}0 & \hspace{0.6cm}0 \\ \hspace{0.6cm}0 & \hspace{0.6cm}0 \end{array}\right]$
        &$\left[ \begin{array}{cc} -0.0605 & -0.0419 \\-0.0517 & -0.0385 \end{array}\right]$
        &$\left[ \begin{array}{cc}  -0.0848 & -0.1243 \\-0.0840 & -0.1399 \end{array}\right]$  \\[3ex]
        &$\tilde{U}(k)$
        &$\left[ \begin{array}{*{2}{p{1.15cm}}} \hspace{0.6cm}0 & \hspace{0.6cm}0 \\ \hspace{0.6cm}0 & \hspace{0.6cm}0 \end{array}\right]$
        &$\left[ \begin{array}{cc} -0.0975 & -0.0525 \\-0.0829 & -0.0630 \end{array}\right]$
        &$\left[ \begin{array}{cc} -0.1902 & -0.1908 \\-0.2225 & -0.2947 \end{array}\right]$  \\[3ex]
        &$V(k)$
        &$\left[ \begin{array}{*{2}{p{1.15cm}}} \hspace{0.6cm}0 & \hspace{0.6cm}0 \\ \hspace{0.6cm}0 & \hspace{0.6cm}0 \end{array}\right]$
        &$\left[ \begin{array}{cc} \ \ 0.0243 & \ \ \ 0.0176 \\  \ \ 0.0298 & \ \ 0.0215 \end{array}\right]$
        &$\left[ \begin{array}{cc} \ \ 0.0090 & \ \ \ 0.0202 \\ \ \ 0.0186 & \ \ \ 0.0422 \end{array}\right]$  \\[3ex]
        &$\tilde{V}(k)$
        &$\left[ \begin{array}{*{2}{p{1.15cm}}} \hspace{0.6cm}0 & \hspace{0.6cm}0 \\ \hspace{0.6cm}0 & \hspace{0.6cm}0 \end{array}\right]$
        &$\left[ \begin{array}{cc} \ \ 0.0905 & \ \ \ 0.0575 \\ \ \ 0.1194 & \ \ \ 0.0769 \end{array}\right]$
        &$\left[ \begin{array}{cc} \ \ 0.0891 & \ \ \ 0.1179 \\ \ \ 0.1550 & \ \ \ 0.2060 \end{array}\right]$  \\[3ex]
        &$P_{1}(k)$
        &$\left[ \begin{array}{cc} -0.0625 & -0.0825 \\-0.0825 & -0.1125 \end{array}\right]$
        &$\left[ \begin{array}{cc} -0.0396 & -0.0593 \\-0.0593 & -0.1129 \end{array}\right]$
        &$\left[ \begin{array}{cc} -0.1141 & -0.1012 \\-0.1012 & -0.1148 \end{array}\right]$  \\[3ex]
        &$Q_{1}(k)$
        &$\left[ \begin{array}{cc} -0.0625 & -0.0825 \\-0.0825 & -0.1125 \end{array}\right]$
        &$\left[ \begin{array}{cc} -0.1286 & -0.1167 \\-0.1167 & -0.1502 \end{array}\right]$
        &$\left[ \begin{array}{cc} -0.3248 & -0.3715 \\-0.3715 & -0.4619 \end{array}\right]$  \\[3ex]
        &$\tilde{P}_{1}(k)$
        &$\left[ \begin{array}{cc} \ \ 1.0625 & \ \ \ 0.0825 \\ \ \  0.0825 & \ \ \  1.1125 \end{array}\right]$
        &$\left[ \begin{array}{cc} \ \ 1.1255 & \ \ \ 0.1195 \\ \ \ 0.1195 & \ \ \ 1.1585 \end{array}\right]$
        &$\left[ \begin{array}{cc} \ \ 1.1729 & \ \ \ 0.2114 \\ \ \ 0.2114 & \ \ \ 1.3676 \end{array}\right]$  \\[3ex]
        &$\tilde{Q}_{1}(k)$
        &$\left[ \begin{array}{cc} \ \ 1.0625 & \ \ \  0.0825 \\ \ \ 0.0825 & \ \  1.1125 \end{array}\right]$
        &$\left[ \begin{array}{cc}\ \ 1.6674 & \ \ \ 0.4629 \\ \ \ 0.4629 & \ \ \ 1.3867 \end{array}\right]$
        &$\left[ \begin{array}{cc} \ \ 2.0130 & \ \ \ 1.2515 \\ \ \ 1.2515 & \ \ \ 2.8022 \end{array}\right]$  \\[3ex]
 \bottomrule
    \end{tabular*}}
\end{center}}
\end{center}
\end{example}

\section{Conclusion}
We have discussed the finite horizon  $H_2/H_\infty$ control problem
of mean-field type for discrete-time systems with state and
disturbance dependent noise. A sufficient condition  has been
derived via the solvability of four coupled matrix-valued equations,
for which, a recursive algorithm has also been provided.   \\

%%Please make sure that your given name is abbreviated as the first capital letter, such as Zhang X T, Tami T,...

\end{document}

%% file: jsscN.tex
\makeatletter
\def\@oddfoot{\hfill}
\newcount\shumeicount
\def\setshumei#1#2#3{%
  \shumeicount=\count0
  \def\@oddhead{%
    \raise-5pt\hbox to0pt{\vrule width\hsize height 0pt depth 0.4pt\hss}\relax
    \ifnum \shumeicount=\count0
      \raise-7pt\hbox to0pt{\vrule width\hsize height 0pt depth 0.4pt\hss}\relax
      #1
    \else
      \ifodd\count0
        #2
      \else
        #3
       \fi
     \fi
  }%
}
\makeatother
\makeatletter
\def\@oddfoot{\hfill}
\newcount\shujiaocount
\def\setshujiao{%
  \shujiaocount=\count0
  \def\@oddfoot{%
      \ifodd\count0
      \else
      \fi
  }%
}
\makeatother
\def\biaoti#1#2#3#4{{
  \vspace*{0.3cm}
  \begin{flushleft} \Large\bf #1\end{flushleft}
  \vspace*{-0.2cm}
      \begin{flushleft}
      \bf #2
      \end{flushleft}
      \footnotetext{\hspace{-6mm} #3\\ #4}}}
\def\dshm#1#2#3#4
{\setshumei{J Syst Sci Complex (#1) #2:
{\thepage--\pageref{LastPage}}\hfill}
            {\hfill {\small #3}\hfill\hbox to0pt{\hss\thepage}}
            {\hbox to0pt{\thepage\hss}\hfill {\small #4}\hfill
            }
            \setshujiao}
\def\drd#1#2
{{\vskip 1cm\small \begin{flushleft}
 #1 \\
 #2\\
\copyright The Editorial Office of JSSC \&  Springer-Verlag Berlin
Heidelberg 2014
\end{flushleft}}}
\def\tilde{\widetilde}
\def\bar{\overline}
\def\epsilon{\varepsilon}
\def\proof{\vspace{1mm}\indent {\it Proof}\quad}